\RequirePackage{fix-cm}
\documentclass[smallextended]{svjour3}       
\smartqed  
\RequirePackage{graphicx}
\RequirePackage{amsfonts}
\RequirePackage{amsmath}
\RequirePackage{amssymb}
\RequirePackage{mathtools}
\RequirePackage{mathabx}
\RequirePackage{csquotes}
\RequirePackage{tipx}
\RequirePackage{textcomp}
\RequirePackage{mathrsfs}
\RequirePackage[margin=1cm]{geometry}

\makeatletter
\newcommand{\eqnum}{\refstepcounter{equation}\textup{\tagform@{\theequation}}}
\makeatother

\begin{document}

\title{An investigation of the non-trivial zeros of the Riemann zeta function}


\author{
Yuri Heymann\        
}


\institute{Yuri Heymann \at
              Georgia Institute of Technology, Atlanta, GA 30332, USA  (Alumni)\\
              \email{y.heymann@yahoo.com}             \\
              \emph{Address in Switzerland:} 3 rue Chandieu, 1202 Geneva. 
 }

\date{}

\maketitle

\begin{abstract}
While many zeros of the Riemann zeta function are located on the critical line $\Re(s)=1/2$, the non-existence of zeros in the remaining part of the critical strip $\Re(s) \in \, ]0, 1[$ is the main scope to be proven for the Riemann hypothesis. The Riemann zeta functional leads to a relation between the zeros on either sides of the critical line. Given $s$ a complex number and $\bar{s}$ its complex conjugate, if $s$ is a zero of the Riemann zeta function in the critical strip $\Re(s) \in \, ]0, 1[$, then $\zeta(s) = \zeta(1-\bar{s})$, as a key proposition to prove the Riemann hypothesis.

\end{abstract}

\section{Introduction}

The Riemann zeta function named after Bernhard Georg Friedrich B. Riemann (1826 - 1866) is an extension of the zeta function to complex numbers, which primary purpose is the study of the distribution of prime numbers \cite{Riemann1859,Mangoldt}. Although the zeta function, an infinite series composed of inverse power semaphores, is part of Euler's product formula establishing the connection between primes and the zeta function, Riemann's approach is connected with more recent developments in complex analysis at the time of Cauchy. The Riemann hypothesis is a statement that all non-trivial zeros lie on the critical line $\Re(s)=1/2$, a foundation of prime-number theory. Not only the Riemann-von Mangoldt explicit formula for the asymptotic expansion of the prime-counting function involves a sum over the non-trivial zeros of the Riemann zeta, but, the Riemann hypothesis has implications for a variety of theorems and the accurate estimate of the error involved in the prime-number theorem more specifically, as seen in \cite{Selberg,SpencerandGraham} (proof and contribution) and earlier work of Hadamard and de la Vall{\'e}e Poussin \cite{Poussin,Hadamard}. The main purpose of aforementioned prime-number theorem is to provide an accurate description of the asymptotic behavior of the prime-counting function $\pi(x)$, representing the number of primes less than or equal to $x$ a variable.

Similar examples of research topics in the realm of prime numbers include the study of functions as $\pi(x)-\pi(x/2) \geq 1,2,3,4,5,..$ when $x \geq 2, 11,17,29,41...$, a step function inspired by the work of Ramanujan et al. (see \cite{Hardy,RamanujanCollPapers}) and prime sequence $sn=2,11,17,...$ made of Ramanujan primes \cite{Sondow}, providing a support for the study of primes, and selected contributions of Hardy and Littlewood, e.g. \cite{HardyLittlewood,Ingham}, etc. 


\vspace{2mm}
\noindent  
With respect to Georg Friedrich prime-number theory, the Riemann zeta function is a holomorphic function defined as follows: 

\begin{equation}
\zeta(s) = \sum_{n=1}^{\infty} \frac{1}{n^s} \,,
\end{equation}

\noindent
where $s$ is a complex number and domain of convergence is $\Re(s) >1$. The domain of convergence of the Riemann zeta function is extended to the left of $\Re(s)=1$ by analytic continuation. The Dirichlet eta function which is the product of the factor $\left(1 - \frac{2}{2^s} \right)$ with the Riemann zeta function is convergent for $\Re(s) > 0$. This alternating series is expressed as follows:

\begin{equation}
\eta(s)= \left( 1 - \frac{2}{2^s}\right) \, \zeta(s) = \sum_{n=1}^{\infty} \frac{(-1)^{n+1}}{n^s} \,,
\end{equation}

\noindent
where $\Re(s)>0$ and $\eta(1)=\ln(2)$ by continuity.

\vspace{2mm}
\noindent
The function $\left( 1- \frac{2}{2^s} \right)$ has an infinity of zeros on the line $\Re(s)=1$ given by $s_k=1+\frac{2 k \pi i}{\ln 2}$ where $k \in \, \mathbb{Z}^{*}$. As $\left( 1 - \frac{2}{2^s} \right) = 2 \times \left( 2^{\alpha - 1} {e}^{i \, \beta \ln 2} - 1\right) / 2^\alpha {e}^{i \, \beta \ln 2}$, the factor $\left( 1 - \frac{2}{2^s} \right)$ has no poles nor zeros in the strip $\Re(s) \in \, ]0,1[$. It follows that by analytic continuation, the Dirichlet eta function can be used as a proxy for zero finding of the Riemann zeta function in the critical strip $\Re(s) \in \, ]0,1[$. The Riemann zetas function can be further extended to the left, i.e. $\Re(s) \leq 0$ by analytic continuation with the Riemann zeta functional. 

\vspace{2mm}
\noindent
In the remaining of the manuscript, the Riemann zeta function is referring to its formal definition and analytic continuation by congruence.

\section{Some elementary propositions}

Before delving into core propositions, a prerequisite is the Riemann zeta functional, which is used explicitly in propositions 2, 4 and 9 and implicitly in 5, 6, 7 and 8. The Riemann zeta functional is expressed as follows:

\begin{equation}
\zeta(s) = \Pi(-s) (2 \pi)^{s-1} 2 \sin \left( \frac{s \pi}{2} \right) \zeta(1-s) \,.
\end{equation}

\noindent
This equation was established by Riemann in 1859. The details of the derivation are provided in \cite{Edwards}, p.13. The standard formulation of (3) is as follows:

\begin{equation}
\zeta(1-s) = \Gamma(s) (2 \pi)^{s-1} 2 \cos \left( \frac{s \pi}{2} \right) \zeta(s) \,,
\end{equation}

\noindent
involves the substitution $s \to 1 - s$ and the relationship $\Pi(s-1)=\Gamma(s)$. Let us proceed with some key propositions.

\vspace{2mm}
\noindent
\textbf{Proposition 1} A formula is introduced here for the calculations further down, which is expressed as follows:

\begin{equation}
\left[ \frac{a^2-b^2}{a^2+b^2} + i \, \frac{-2 a b}{a^2+b^2}  \right] \left[a + i \, b\right] = a - i \, b \,.
\end{equation}

\noindent
where $a$ and $b$ are reals.

\begin{proof}
By identifying the real and imaginary parts of $(x + i \, y) (a + b \, i)= a - i \, b$, we obtain two equations $a x - b y = a$ and $b x + a y = -b$. Eq. (5) stems from these two equations.
\end{proof}

\vspace{2mm}
\noindent
\textbf{Proposition 2} Given $s$ a complex number and $\bar{s}$ its complex conjugate, we have:

\begin{equation}
\zeta(\bar{s}) = \bar{\zeta}(s) \,,
\end{equation}

\noindent
where $\bar{\zeta}(s)$ is the complex conjugate of $\zeta(s)$.

\begin{proof}
Let us say $s= \alpha + i \, \beta$ where $\alpha$ and $\beta$ are real numbers and $i^2 = -1$. We have:

\begin{equation}
\zeta(s) = \sum_{n=1}^{\infty} \frac{[\cos(\beta \, \ln n) - i \, \sin(\beta \, \ln n)]}{n^\alpha} \,,
\end{equation}

\noindent
where $\Re(s) >1$.

\noindent
We have $\frac{1}{n^s}= \frac{1}{n^\alpha \, {e}^{i \, \beta \ln n}}$. We then multiply both the numerator and denominator by $\cos(\beta \ln n) - i \, \sin(\beta \ln n)$. After several simplifications, we get (7). Note that in (7) when $\beta$ changes its sign, the real part of $\zeta(s)$ remains unchanged, while imaginary part changes its sign. This means, the Riemann zeta function has mirror symmetry with respect to the real axis of the complex plane and $\zeta(\bar{s})= \bar{\zeta}(s)$ when $\Re(s) > 1$.

\vspace{2mm}
\noindent
The same way, the Dirichlet eta function is rewritten as follows:

\begin{equation}
\eta(s) = \sum_{n=1}^{\infty} (-1)^{n+1} \frac{[\cos(\beta \ln n) - i \, \sin(\beta \ln n)]}{n^\alpha} \,,
\end{equation}

\noindent
where $\Re(s) >0$ and $\eta(1)= \ln 2$ by continuity.

\vspace{2mm}
\noindent
Using the same reasoning as above, the Dirichlet eta function is symmetrical with respect to the real axis of the complex plane. Furthermore, the factor $(1-2^{1-s}) = 1 - 2^{1 - \alpha} (\cos(\beta \ln 2) - i \, \sin(\beta \ln 2))$ has mirror symmetry with respect to the real axis. As the product of holomorphic functions symmetrical with respect to the real axis yields a similar holomorphic function, also symmetrical with respect to the real axis, we get $\zeta(\bar{s})=\bar{\zeta}(s)$ when $\Re(s) \in \, ]0,1[$.

\vspace{2mm}
\noindent
Let us consider the Riemann zeta functional (4) and introduce $\xi : \mathbb{C} \to \mathbb{C}$ defined as:

\begin{equation}
\xi(s) = 2 \Gamma(s) (2 \pi)^{-s} \, \cos(\frac{s \pi}{2}) \,.
\end{equation}

\noindent
For the Riemann zeta function to be symmetrical with respect to the real axis when $\Re(s) < 0$, by analytic continuation aforementioned $\xi(s)$ function has to be symmetrical as well.

From formula 6.1.23 p.256 in \cite{Abramowitz}, we have $\bar{\Gamma}(s) = \Gamma(\bar{s})$. Furthermore  $(2\pi)^{-s}$ is symmetrical with respect to the real axis. From 4.3.56 p.74 in \cite{Abramowitz}, $\cos \left( \frac{s \pi}{2} \right) = \cos \left( \frac{\alpha \pi}{2} \right) \cosh \left( \frac{\beta \pi}{2}\right) - i \, \sin \left( \frac{\alpha \pi}{2} \right) \sinh\left( \frac{\beta \pi}{2} \right)$. As $\cosh(-x)=\cosh(x)$ and $\sinh(-x) = -\sinh(x)$, $\cos \left( \frac{s \pi}{2}\right)$ is symmetrical with respect to the real axis. Therefore, $\zeta(\bar{s})=\bar{\zeta}(s)$ when $\Re(s) < 0$. As the Riemann zeta function is a meromorphic function, by continuity $\zeta(\bar{s})=\bar{\zeta}(s)$ on the lines $\Re(s)=0$ and $\Re(s)=1$.
\end{proof}

\vspace{2mm}
\noindent
\textbf{Proposition 3} Given $s$ a complex number and $\bar{s}$ its complex conjugate, we define a holomorphic function $\nu : \mathbb{C} \to \mathbb{C}$ such that:

\begin{equation}
\zeta(s) = \nu(s) \zeta(1 - \bar{s}) \,,
\end{equation}

\noindent
where $\nu(s)$ is equal to one when $\Re(s)=1/2$.

\vspace{2mm}
\noindent
\textit{Rationale:}
Let us set $s = \alpha + i \, \beta$ where $\alpha$ and $\beta$ are reals. When $\alpha = 1/2$, we have $\zeta(s)=\zeta(1/2 + i \, \beta)$ and $\zeta(1- \bar{s})= \zeta(1/2 + i \, \beta)$, meaning $\zeta(s)=\zeta(1- \bar{s})$ when $\Re(s)=1/2$. For a given $\alpha \neq 1/2$, we introduce a complex factor $\nu(s)$ in (10) because the two equalities $\Re(\zeta(s))=\Re(\zeta(1 - \bar{s}))$ and $\Im(\zeta(s))=\Im(\zeta(1 - \bar{s}))$ are no longer satisfied with one degree of freedom alone given by $\beta$.

\vspace{2mm}
\noindent
\textbf{Proposition 4} Given $s$ a complex number and $\bar{s}$ its complex conjugate, we have:

\begin{equation}
\frac{1}{\nu(s)}=2 \, \frac{\Gamma(\bar{s})}{(2 \pi)^{\bar{s}}} \, \cos\left( \frac{\pi \bar{s}}{2} \right) \, \left[ \frac{u^2 - v^2}{u^2+v^2} + i \, \left( \frac{-2 u v}{u^2+v^2} \right) \right] \,,
\end{equation}

\noindent
where $u=\Re(\zeta(s))$ and $v=\Im(\zeta(s))$.

\begin{proof}
This formula stems from the Riemann zeta functional. From \textit{proposition no 3}, we have $\frac{1}{\nu(s)}=\frac{\zeta(1 - \bar{s})}{\zeta(s)}$. From the Riemann zeta functional (4) expressed as $\zeta(1-s)=\frac{\Gamma(s)}{(2 \pi)^s} 2 \cos \left( \frac{\pi s}{2} \right) \zeta(s)$, we have $\zeta(1- \bar{s})=\frac{\Gamma(\bar{s})}{(2 \pi)^{\bar{s}}} 2 \cos \left( \frac{\pi \bar{s}}{2}\right) \zeta(\bar{s})$. From \textit{proposition no 2}, $\zeta(\bar{s})=\bar{\zeta}(s)$, where $\bar{\zeta}(s)$ is the complex conjugate of $\zeta(s)$.\\
By \textit{proposition no 1} we get $\zeta(\bar{s})=\zeta(s) \left[ \frac{u^2 - v^2}{u^2+v^2} + i \, \left( \frac{-2 u v}{u^2 + v^2} \right) \right]$ where $u=\Re(\zeta(s))$ and $v=\Im(\zeta(s))$. Eq. (11) follows.
\end{proof}

\vspace{2mm}
\noindent
\textbf{Proposition 5} Given $s$ a complex number, the function $\nu(s)$ as defined above is equal to zero only at points $s=0, -2, -4,-6,-8,...,-n$ where $n$ is an even integer.

\begin{proof}
We set $s= \alpha + i \, \beta$ where $\alpha$ and $\beta$ are reals. The function $\nu(s)$ tends to zero, when its reciprocal $\frac{1}{\nu(s)}$ approaches $\pm \infty$. We invoke \textit{proposition no 4} to find the values where $\frac{1}{\nu(s)} \to \pm \infty$.

\vspace{2mm}
\noindent
(i) We note that the term $\left[ \frac{u^2-v^2}{u^2+v^2} + i \, \left( \frac{-2 u v}{u^2 + v^2} \right) \right]$ is bounded: we have $-1 \leq \frac{-v^2}{u^2+v^2} \leq \frac{u^2-v^2}{u^2+v^2} \leq \frac{u^2}{u^2+v^2} \leq 1$ and $-2 \leq \frac{-2 \text{max}(|u|,|v|)^2}{u^2+v^2} \leq \frac{-2 u v}{u^2+v^2} \leq \frac{2 \text{max}(|u|,|v|)^2}{u^2+v^2} \leq 2$.

\vspace{2mm}
\noindent
(ii) We note that the expression $\left[ \frac{u^2-v^2}{u^2+v^2} + i \, \left( \frac{-2 u v}{u^2 + v^2} \right)\right]$ is never equal to zero as it is not possible to have both the real and imaginary parts equal to zero at the same time. For the real part to be equal to zero either $u=v$ or $u=-v$. When $u=v$ (or $u = -v$), the imaginary part is equal to $-1$ (or $1$). For the imaginary part to be equal to zero, either $u$ or $v$ should be equal to zero. If $u$ is equal zero (or $v$ is equal to zero), then the real part is equal is equal to $-1$ (or $1$).

\vspace{2mm}
\noindent
Eq. (5) can also be expressed as $\left[\frac{u^2-v^2}{u^2+v^2} + i \, \left( \frac{-2u v}{u^2+v^2} \right) \right]=\frac{u - i v}{u + i v}$. This formula is properly defined when $u+i v \neq 0$.

\vspace{2mm}
\noindent
As implied by (i) and (ii), when $u + i \, v$ is in the neighborhood of 0, no matter how close to 0, corresponding factor $\frac{u - i v}{u + i v}$ cannot take the value zero and is bounded.

\vspace{2mm}
\noindent
Consequently of (i) and (ii), $\nu(s)=0$ at a point $s_0$ if and only if $\lim_{s \to s_0} \, \cos \left( \frac{\pi \bar{s}}{2}\right) \, \Gamma(\bar{s})=\pm \infty$.

\vspace{2mm}
\noindent
We check that $\cos \left( \frac{\pi \bar{s}}{2}\right)$ is bounded. We have $\cos \left( \frac{\pi}{2} (\alpha - i \, \beta)\right) = \cos \left( \frac{\pi \alpha}{2}\right) \cosh \left( \frac{-\pi \beta}{2} \right) - i \, \sin \left( \frac{\pi \alpha}{2} \right) \, \sinh \left( \frac{-\pi \beta}{2}\right)$. As the cosine, the sine, the hyperbolic cosine and the hyperbolic sine are bounded on $]-\infty, \infty[$, $\cos \left( \frac{\pi \bar{s}}{2}\right)$ is bounded when $\beta$ is finite.

\vspace{2mm}
\noindent
Hence, the condition $\nu(s)=0$ occurs only if $\Gamma(\bar{s}) \to \pm \infty$ or when the reciprocal Gamma $\frac{1}{\Gamma(\bar{s})}=0$. The Euler form of the reciprocal Gamma function (see \cite{Abramowitz} p.8) is $\frac{1}{\Gamma(s)}= s \, \prod_{n=1}^{\infty} \frac{1+\frac{s}{n}}{\left( 1 + \frac{1}{n} \right)^s}$. Thus, the full set of zeros of the reciprocal Gamma is given by $s=0, -1, -2, -3,...,-n$ where $n \in \, \mathbb{N}$. These are the points where $\nu(s)$ can potentially be equal to zero.

\vspace{2mm}
\noindent
We check the points where the cosine term $\cos \left( \frac{\pi \bar{s}}{2} \right)$ is equal to zero. We have $\cos \left( \frac{\pi \bar{s}}{2}\right)= \cos \left(\frac{\pi}{2} \left( \alpha - i \, \beta\right)\right)= \cos \left( \frac{\pi \alpha}{2}\right) \, \cosh \left( \frac{-\pi \beta}{2}\right) - i \, \sin \left( \frac{\pi \alpha}{2} \right) \, \sinh \left( \frac{- \pi \beta}{2} \right)$. As the hyperbolic cosine is never equal to zero, the real part of $\cos \left( \frac{\pi \bar{s}}{2} \right)$ is equal to zero if and only if $\cos \left( \frac{\pi \alpha}{2} \right)=0$. As we cannot have both the sine and cosine functions equal to zero simultaneously, when they share the same argument, the imaginary part of $\cos \left( \frac{\pi \bar{s}}{2} \right)$ equals zero if and only if $\sinh \left(\frac{-\pi \beta}{2} \right)=0$, which occurs if and only if $\beta=0$. Hence, the term $\cos \left( \frac{\pi \bar{s}}{2} \right)$ can be equal to zero on the real line at the points $s=\pm 1, \pm 3, \pm 5, \pm 7, ...$ Thus, $\nu(s)$ is equal to zero with certainty at the points $s=0, -2, -4, -6, -8,...$

\vspace{2mm}
\noindent
What about the points $s = -1, -3, -5, -7, -9, ...$ ?

\vspace{2mm}
\noindent
The Euler form of the Gamma function is expressed as $\Gamma(s)= \frac{1}{s} \, \prod_{n=1}^{\infty} \, \frac{\left(1 + \frac{1}{n} \right)^s}{1 + \frac{s}{n}}$. The Taylor series of $\cos \left( \frac{\pi s}{2}\right)$ in $-m$ where $m$ is an odd integer is as follows: $\cos \left( \frac{\pi s}{s}\right)=\frac{\pi}{2} \, \sin \left( \frac{m \pi}{2} \right) \left( s + m\right) - \left( \frac{\pi}{2} \right)^3 \, \sin \left( \frac{m \pi}{2} \right) \, \left( s + m\right)^3+...$ If we multiply the first term of the Taylor series with the m-th product of the Euler form of the Gamma function, we obtain a factor of order $\frac{s+m}{s+m}$, which is finite when $s$ tends to $-m$. The multiplication of the successive terms of the Taylor series with the m-th product of the Euler form of the Gamma function leads to a factor of order $\frac{(s+m)^k}{s+m}$ where $k=3,5,7,...,$ which converges towards zero when $s$ tends to $-m$. Furthermore, the limit of the n-th product of the Euler's form of the Gamma function $\frac{\left( 1 + \frac{1}{n}\right)^s}{1 + \frac{s}{n}}$ when $n$ tends to $+\infty$ is equal to $1$. Therefore, $\lim_{s \to s_0} \, \cos \left( \frac{\pi \bar{s}}{2} \right) \, \Gamma(\bar{s})$ is finite at the points $s= -1, -3, -5, -7, -9, ...$ so we can say that $\nu(s) \neq 0$ at these points.

\vspace{2mm}
\noindent
As  a validation step, from $\textit{proposition no 3}$ we have $\zeta(s)=\nu(s) \zeta(1-\bar{s})$, and from the proposition that the only pole of $\zeta(z)$ is at the point $z=1$. Furthermore, according to §3, we have $\zeta(1-\bar{s}) \neq 0$ when $\Re(s) < 0$. Hence, $\nu(s)$ cannot be equal to zero at the points $s=-1,-3,-5,-7,-9,...$; otherwise, $\zeta(s)$ would be equal to zero at these points, which contradicts the non-trivial zeros of the Riemann zeta function (see §4). Although $\nu(0)=0, \zeta(0) \neq 0$ because the function $\zeta(z)$ has a pole at $z=1$ (see \textit{prop. no. 9}).

\end{proof}

\vspace{2mm}
\noindent
\textbf{Proposition 6} Given $s$ a complex number, the function $\nu(s)$ defined above has an infinite number of poles at the points $s=1, 3, 5, 7, 9,..., n$ where $n$ is an odd integer.

\begin{proof}
We use the expression of the reciprocal of $\nu(s)$ introduced in \textit{proposition no 4}. We have shown in proof of \textit{proposition no 5}, that the factor $\left[ \frac{u^2-v^2}{u^2+v^2} + i \, \left( \frac{-2 u v}{u^2+v^2}\right)\right]$ cannot take the value zero and is bounded. Hence, a point $s_0$ is a pole of the function $\nu(s)$ if and only if $\lim_{s \to s_0} \, \Gamma(\bar{s})\, \cos \left( \frac{\pi \bar{s}}{2} \right)=0$. Such poles occur either when $\Gamma(\bar{s})$ or $\cos \left( \frac{\pi \bar{s}}{2} \right)$ are equal to zero. The term $\cos \left( \frac{\pi \bar{s}}{2}\right)$ is equal to zero at the points $s=\pm 1, \pm 3, \pm 5, \pm 7,...$ The Euler form of the Gamma function is expressed as $\Gamma(s)=\frac{1}{s} \, \prod_{n=1}^{\infty} \frac{\left( \right)^s}{1+\frac{s}{n}}$. As $\forall n \in \, \mathbb{N}^{*}$ the term $1+ \frac{1}{n}$ is not equal to zero, we can say that the Gamma function is never equal to zero. However, the Gamma function tends to infinity at the points $s=0, -1, -2, -3, ...,-n$ where $n \in \mathbb{N}$, hence the points $s=1,3,5,7,...$ can be determined to be poles of $\nu(s)$ with certainty. The points $s=-1,-3,-5,-7,...$ are poles of $\nu(s)$ only if the limit of $\Gamma(s) \, \cos \left( \frac{\pi s}{2}\right)$ is equal to zero when $s$ tends to either of these points. The Taylor series of $\cos \left( \frac{\pi s}{2}\right)$ in $-m$ where $m$ is an odd integer is as follows: $\cos \left( \frac{\pi s}{2}\right)=\frac{\pi}{2} \sin \left( \frac{\pi m}{2} \right) \left( s + m \right) - \left( \frac{\pi}{2} \right)^3 \sin \left( \frac{\pi m}{2} \right) \left( s + m\right)^3+...$ If we multiply the m-th product of the Euler's form of the Gamma function with the first term of the Taylor series we computed for $\cos \left( \frac{\pi s}{2} \right)$, we obtain a factor of order $\frac{s+m}{s+m}$, which does not converge towards zero when $s$ tends to $-m$. In addition, the limit of the n-th product of the Euler's form of the Gamma function $\frac{\left( 1 + \frac{1}{n} \right)^s}{1+\frac{s}{n}}$ when $n$ tends to $+\infty$ is equal to $1$. Hence, $\lim_{s \to s_0} \, \Gamma(s) \cos \left( \frac{\pi s}{2} \right) \neq 0$ at the points $s_0 = -1,-3,-5,-7,...$ Thus, the points $s=-1,-3,-5,-7,...$ are not within the set of poles of $\nu(s)$.
\end{proof}

\vspace{2mm}
\noindent
\textbf{Proposition 7} Given $s$ a complex number and $\bar{s}$ its complex conjugate, if $s$ is a complex zero of the Riemann zeta function in the strip $\Re{s} \in \, ]0,1[$, then $1-\bar{s}$ must be a zero.

\begin{proof}
From $\textit{proposition 5}$, $\nu(s)$ is only equal to zero at even negative integers including zero, and never reaches zero when $\Re(s) > 0$. Furthermore, in $\textit{proposition no 6}$ we have shown that $\nu(s)$ has no poles in the strip $\Re(s) \in \, ]0,1[$. As $\zeta(s)=\nu(s) \zeta(1-\bar{s})$ from $\textit{proposition no 4}$ and given that $\nu(s)$ has no poles nor zeros in the strip $\Re(s) \in ]0,1[$ (see \textit{prop. no. 6}), it follows that if $s$ is zero, then $1-\bar{s}$ must also be a zero.
\end{proof}

\vspace{2mm}
\noindent
\textbf{Proposition 8} Given $s$ a complex number and $\bar{s}$ its complex conjugate, if $s$ is a complex zero of the Riemann zeta function in the strip $\Re{s} \in \, ]0,1[$, then we have:

\begin{equation}
\zeta(s)=\zeta(1 - \bar{s}) \,.
\end{equation} 

\begin{proof}
As a corollary of $\textit{proposition no 7}$, if $s$ and $1-\bar{s}$ are zeros of the Rieman zeta function, both $\zeta(s)$ and $\zeta(1-\bar{s})$ are equal to zero meaning both terms are matching. The converse is not necessarily true.  
\end{proof}

\noindent
\textit{Remark:}
Non-transitivity of the Riemann zeta fonctional is tightly linked with the statement that there is a finite number of roots in the critical strip not on $\Re(s)=1/2$, such that $\left| \sum_{n=2}^{\infty} (-1)^{n+1} \frac{e^{i\, \beta \ln n}}{n^{\alpha}}\right| = \left| \sum_{n=2}^{\infty} (-1)^{n+1} \frac{e^{i \, \beta \ln n}}{n^{1-\alpha}}\right|=1$, and that each of there roots is not a zero of the Riemann zeta function. We suppose there exists such a root in the neighborhood of the conjugated pair of point $s_0=0.7332 + i \, 289.7999$ and $s_1=0.2668 + i \, 289.7999$.

\section{A recap on the zeros of the Riemann zeta function}

\subsection{Poles of the Riemann zeta function}

\textbf{Proposition 9}
The only pole of the Riemann zeta function $\zeta(s)$ is a simple pole at $s=1$.

\vspace{2mm}
\noindent
This is a theorem documented in \cite{Titchmarsh} p.13, see the below sketch.  

\noindent
Using Abel's lemma for summation by parts, see \cite{Krantz} at p.58, we get:

\begin{equation}
\sum_{n=1}^{m} n^{-s} = \sum_{n=1}^{m-1} n \, \left(n^{-s} - (n+1)^{-s} \right) + m^{1-s} \,,
\end{equation}

\noindent
where $m \in \, \mathbb{N}$.

\vspace{2mm}
\noindent
When $\Re(s)>1$, we have $\lim\limits_{m \to \infty} m^{1-s}=0$. Thus, we get:

\begin{equation}
\begin{split}
\sum_{n=1}^{\infty} n^{-s} & = \sum_{n=1}^{\infty} n \left( n^{-s} -(n+1)^{-s} \right) \\
& = s \sum_{n=1}^{\infty} n \int_{n}^{n+1} x^{-s-1} \, dx \\
& = s \int_{1}^{\infty} \lfloor x \rfloor \, dx \\
& = \frac{s}{s-1} - s \int_{1}^{\infty} \{x\} x^{-s-1} \, dx \,,
\end{split}
\end{equation}

\noindent
where $\lfloor x \rfloor$ denotes the floor function of $x$ and $\{x\}=x-\lfloor x \rfloor$ the fractional part of $x$.

\vspace{2mm}
\noindent
Hence, we have:

\begin{equation}
\zeta(s)=\frac{s}{s-1} - s \int_{1}^{\infty} \{ x \} x^{-s-1} \, dx \,,
\end{equation}

\noindent
where $\Re(s) > 1$.

\vspace{2mm}
\noindent
In order for $s=1$ to be the only pole of $\zeta(s)$ on the half-plane $\Re(s) \geq 1$, we have to show that the integral $\int_{1}^{\infty} \{ x\} x^{-s-1} \, dx$ in (15) is bounded when $\Re(s) \geq 1$.

\vspace{2mm}
\noindent
Given a sequence $u_i$ of complex numbers, lower triangle inequality yields $|\sum_i u_i| \leq \sum_i |u_i|$. As an integral can be represented as an infinite sum,  this inequality still holds and we get:

\begin{equation}
\begin{split}
\left| \int_{1}^{\infty} \{ x \} x^{-s-1} \, dx \right| & \leq \int_{1}^{\infty} \left| \{ x\} x^{-s-1}\right| \, dx \\
& \int_{1}^{\infty} \left| \{ x \} x^{(-\alpha -1)} {e}^{-i \, \beta \ln x}\right| \, dx \\
&  \int_{1}^{\infty} \{ x \} x^{-\alpha -1} \, dx
\end{split}
\end{equation}

\noindent
As $\int_{n}^{n+1} \{ x \} x^{-\alpha -1} < \int_{n}^{n+1} x^{-\alpha -1} \, dx$, we get:

\begin{equation}
\begin{split}
\left|\int_{1}^{\infty} \{ x \} x^{-s-1} \, dx \right| & < \sum_{n=1}^{\infty} \int_{n=1}^{n+1} x^{-\alpha-1} \, dx \\
& < \frac{1}{\alpha} \sum_{n=1}^{\infty} \left( n^{-\alpha} - (n+1)^{-\alpha}\right) \,.
\end{split}
\end{equation}

\noindent
When $\alpha=1$, we have:

\begin{equation}
\begin{split}
\left| \int_{1}^{\infty} \{ x \} x^{-s-1} \, dx\right| & < \sum_{n=1}^{\infty} \left( n^{-1} - (n+1)^{-1} \right) \\
& < \sum_{n=1}^{\infty} \frac{1}{n (n+1)} \,.
\end{split}
\end{equation}

\noindent
From the integral test (\cite{Grigorieva} p. 132), the series $\sum_{n=1}^{\infty} \frac{1}{n(n+1)} < \sum_{n=1}^{\infty} \frac{1}{n^2}$ is convergent. Thus, the only pole on the line $\Re(s) =1$ is at $s=1$, which appears to be a simple pole.

\vspace{2mm}
\noindent
From integral test, when $\alpha > 1$, both series $\sum_{n=1}^{\infty} \frac{1}{n^{\alpha}}$ and $\sum_{n=1}^{\infty} \frac{1}{(n+1)^{\alpha}}$ are convergent. Hence, $\left| \int_{1}^{\infty} \{ x \} x^{-s-1} \, dx\right|$ in (18) is bounded. Thus, we can say that $\zeta(s)$ has no poles when $\Re(s) > 1$.

%

\vspace{2mm}
\noindent
To show there are no poles for negative $\alpha$, we use the property of the Gamma function: 

\begin{equation}
\Gamma(s)=\frac{1}{s} \prod_{n=1}^{\infty} \frac{\left( 1 + \frac{1}{n}\right)^s}{1 + \frac{s}{n}} \,,
\end{equation}

\noindent
that $\Gamma$ has no poles when $\Re(s) \geq 1$ and cosine term in (4) is bounded. Also, we have shown that the only pole of $\zeta(s)$ in the half-plane $\Re(s) \geq 1$ is a point at $s=1$. The zeta function evaluated in $0$ is $\zeta(0)=-\frac{1}{2}$ see \cite{Veen} p. 135. Thus, by extrapolation from the Riemann zeta functional (4) when $\alpha$ positive, we can say that $\zeta(s)$ has no poles when $\Re(s) \leq 0$. 

\vspace{2mm}
\noindent
Using the Dirichlet eta function as an extension of the Riemann zeta function, given the factor $\left( 1 - \frac{2}{2^{s}}\right)$ has no poles nor zeros in the strip $\Re(s) \in \, ]0,1[$, one has to show that the Dirichlet eta function is bounded, to say the $\zeta(s)$ has no poles when $\Re(s) \in \, ]0,1[$. The Dirichlet eta function in its integral form is expressed as follows:

\begin{equation}
\eta(s) = \frac{1}{\Gamma(s)} \, \int_{0}^{\infty} \frac{x^{s-1}}{{e}^{x}+1} \, dx \,.
\end{equation}

\vspace{2mm}
\noindent
As the Gamma function is never equal to zero (see in prop. no 6), $\eta(s)$ is bounded if the integral $\int_{0}^{\infty} \frac{x^{s-1}}{{e}^{x}+1} \, dx$ is bounded. Let us consider the case when $\alpha \in \, ]0,1[$. Since:

\begin{equation}
\begin{split}
\left| \int_{0}^{\infty} \frac{x^{s-1}}{{e}^{x}+1} \right| & \leq \int_{0}^{\infty} \left| \frac{x^{s-1}}{{e}^{x}+1} \right| \, dx \\
& \leq \int_{0}^{\infty} \frac{x^{\alpha-1}}{{e}^{x}+1} \, dx \\
& < \int_{0}^{1} \frac{x^{\alpha -1}}{{e}^{x}+1} \, dx + \int_{1}^{\infty} \, \frac{1}{{e}^{x}+1} \, dx \\
& < \int_{0}^{1} \frac{x^{\alpha-1}}{{e}^{x}+1} \, dx + \int_{1}^{\infty} \, {e}^{-x} \, dx
\end{split}
\end{equation}

\noindent
We solve $\int_{0}^{1} \frac{x^{\alpha-1}}{{e}^{x}+1} \, dx$ using integration by part $\int \, u^{\prime} \, v = [u \, v] - \int \, u v^{\prime}$ with $u^{\prime}=x^{\alpha -1}$ and $v = \frac{1}{{e}^{x}+1}$. We get:

\begin{equation}
\begin{split}
\int_{0}^{1} \frac{x^{\alpha-1}}{{e}^{x}+1} \, dx & = \left[ \frac{x^{\alpha}}{\alpha} \frac{1}{{e}^{x}+1} \right]_{0}^{1} + \int_{0}^{1} \frac{x^{\alpha}}{\alpha} \, \frac{{e}^{x}}{({e}^{x}+1)^2} \, dx \\
& < \left[ \frac{x^{\alpha}}{\alpha} \frac{1}{{e}^{x}+1} \right] _{0}^{1}+ \int_{0}^{1} \frac{1}{\alpha} \frac{{e}^x}{({e}^x+1)^2} \, dx \\
& < \frac{1}{\alpha \left( e^1 +1 \right)} + \frac{1}{\alpha} \left( \frac{1}{2} - \frac{1}{{e}^1 +1} \right) \,.
\end{split}
\end{equation}

\noindent
Thus, we have:

\begin{equation}
\left| \int_{0}^{\infty} \frac{x^{s-1}}{{e}^{x}+1} \, dx \right| < \frac{1}{\alpha \left({e}^1 + 1 \right)} + \frac{1}{\alpha} \left(\frac{1}{2} -\frac{1}{{e}^1+1} \right) + {e}^1 \,,
\end{equation}

\noindent
and $\eta(s)$ is bounded in the strip $]0,1[$, so we can say that $\zeta(s)$ has no poles in the critical strip.

\subsection{Zeros when $\Re(s) > 1$}

\noindent
A sketch of the method for the zeros when $\Re(s) >1$ (see \cite{Titchmarsh}) is made available below. The method is based on Euler product and the integral form of the remainder of the Taylor expansion of ${e}^{x}$.

\vspace{2mm}
\noindent
The Taylor expansion of ${e}^{x}$ in $0$ may be expressed as follows:

\begin{equation}
{e}^{x} = 1 + x + \varepsilon(x) \,,
\end{equation}

\noindent
where $\varepsilon(x)$ is the remainder of the Taylor approximation.

\vspace{2mm}
\noindent
The integral form of the remainder of the Taylor expansion of ${e}^{x}$ is as follows:

\begin{equation}
\varepsilon(x) = \int_{0}^{x} \left( x - u \right) \, {e}^u \, du \,.
\end{equation}

\noindent
When $u \in \, [0,x]$, we have $(x - u) {e}^{u} \geq 0$. Hence, we have $\varepsilon(x) \geq 0$.  We get:

\begin{equation}
{e}^{x} \geq 1 + x \,.
\end{equation}

\noindent
Given $p$ a prime number larger than one, and $s=\alpha +i \, \beta$ a complex number where $\alpha$ and $\beta$ are reals. We have:

\begin{equation}
\begin{split}
\left| 1 - p^{-s} \right| & \leq 1 + \left| p^{-s}\right| \\
& \leq 1 + \left| p^{-\alpha -i \, \beta}\right| \\
& \leq 1 + p^{-\alpha} \times \left|{e}^{(i \, \beta \ln p)} \right| \\
& \leq 1 + p^{-\alpha} \,.
\end{split}
\end{equation}

\noindent
By setting $x = p^{-\alpha}$ in (26), we get:

\begin{equation}
1 + p^{-\alpha} \leq \text{exp} \left( p^{-\alpha}\right) \,.
\end{equation}

\noindent
From (27) and (28), we get:

\begin{equation}
\left| 1 - p^{-s} \right| \leq \text{exp} \left( p^{-\alpha}\right) \,.
\end{equation}

\noindent
The Euler product \cite{Edwards} p.6 is expressed as:

\begin{equation}
\zeta(s) = \prod_{p}^{\infty} \frac{1}{1- p^{-s}} \,,
\end{equation}

\noindent
where $p$ is the sequence of prime numbers larger than $1$ and $\Re(s) > 1$.

\vspace{2mm}
\noindent
Hence, we have:

\begin{equation}
\left| \zeta(s) \right| = \prod_{p}^{\infty} \frac{1}{\left| 1 - p^{-s}\right|} \,.
\end{equation}

\noindent
Using inequality (29) in (31), we get:

\begin{equation}
\begin{split}
\left| \zeta(s)\right| & \geq \prod_{p}^{\infty} \text{exp} \left( -p^{-\alpha}\right) \\
& \geq \text{exp} \left( - \sum_{p}^{\infty} p^{-\alpha} \right) \\
& > 0.  
\end{split}
\end{equation}

\noindent
We have $0 < \sum_{p}^{\infty} \, p^{-\alpha} < \sum_{n=1}^{\infty} \, n^{-\alpha}$. From integral test, the series $\sum_{n=1}^{\infty} n^{-\alpha}$ is convergent when $\alpha > 1$, hence $\text{exp} \left(-\sum_{p}^{\infty} \right) > 0$. From (32), we get $\forall \alpha > 1$, $\left|\zeta(s) \right| > 0$. Hence, we can say the Riemann zeta function has no zeros when $\Re(s) >1$.

\subsection{Zeros when $\Re(s) <0$}

\noindent
The below functional equation is used for the zeros of the Riemann zeta function when $\Re(s) <0$:

\begin{equation}
\pi^{-s/2} \Gamma\left( \frac{s}{2}\right) \zeta(s) = \pi^{-(1-s)/2} \Gamma\left( \frac{1-s}{2}\right) \, \zeta \left( 1 - s \right) \,. 
\end{equation}

\noindent
The derivation of the above variant of the Riemann zeta functional is provided in \cite{Ivic}, p. 8-12. The Gamma function is never equal to zero (see insider proof \textit{prop. no 6}), hence the zeros of the Riemann zeta function when $\Re(s) < 0$ are found at the poles of $\Gamma\left( \frac{s}{2}\right)$ provided $\zeta(1-s)$ has no poles when $\Re(s) < 0$ (see \textit{prop. no. 9}). 

\vspace{2mm}
\noindent
Because the Riemann zeta function has no zeros when $\Re(s) >1$ (see §3.2), when $\Re(s) >1$, $\zeta(1-s)$ in (33) is never equal to zero. Moreover, the Gamma function is never equal to zero (see insider proof prop. no 6), meaning that the zeros of the Riemann zeta function when $\Re(s) <0$ are found at the poles of $\Gamma \left( \frac{s}{2} \right)$ provided $\zeta(1-s)$ is not a pole (which is true from §3.2 when $s \neq 0$). Thus, we obtain the  trivial zeros of $\zeta(s)$ at $s = -2, -4, -6, -8,...$ (see \textit{prop. no 5}). The point $s=0$ is not a zero due to the pole at $\zeta(1)$. We have $\zeta(0)= -\frac{1}{2}$, see \cite{Veen}.

\subsection{Zeros on the line $\Re(s)=0,1$}

The elaboration that Riemann zeta function has no zeros on the line $\Re(s) = 1$ is a result of Hadamard and de la Vall{\'e}e Poussin respective proofs of the prime-number theorem \cite{Poussin,Hadamard}. A sketch of the proof of de la Vall{\'e}e Poussin is available below.
\vspace{2mm}
\noindent
From Euler's product, we have:

\begin{equation}
\frac{1}{\zeta(s)} = \prod_{p}^{\infty} \left( 1 - p^{-s}\right) \,
\end{equation}

\noindent
where $p$ is the sequence of prime numbers larger than $1$. Eq. (34) is defined for $\Re(s) > 1$.

\vspace{2mm}
\noindent
By expanding (34), we have:

\begin{equation}
\frac{1}{\zeta(s)}= 1 - \sum_{p}^{\infty} \left( p^{-s}\right) + \sum_{p < q}^{\infty} \left( p q\right)^{-s} - \sum_{p < q < r}^{\infty} \left( p q r\right)^{-s}+... \,,
\end{equation}

\noindent
where $p, q, r, ...$ are prime numbers. Hence, we get an infinite sum of integers, composed of the product of unique primes. For such integer $n$, the coefficient of $n^{-s}$ is $+1$ if the number of prime factors of $n$ is even, and $-1$ if odd. Thus, we get:

\begin{equation}
\frac{1}{\zeta(s)} = \sum_{n=1}^{\infty} \mu(n) \, n^{-s} \,,
\end{equation}

\noindent
where $\mu(n)$ is the M{\"o}bius function.

\vspace{2mm}
\noindent
We have $\zeta(s) = \sum_{n=1}^{\infty} \, \frac{1}{n^{s}} = \sum_{n=1}^{\infty} \, {e}^{-s \, \ln (n)}$. Hence, we get:

\begin{equation}
\zeta^{\prime}(s) = - \sum_{n=1}^{\infty} \, \ln (n) \, n^{-s} \,.
\end{equation}

\noindent
The product of (36) with (37) leads to the logarithmic derivative of $\zeta(s)$, defined for $\Re(s) > 1$. The expansion yields:

\begin{equation}
\frac{\zeta^{\prime}(s)}{\zeta(s)} = - \sum_{n=2}^{\infty} \Lambda(n) \, n^{-s} \,,
\end{equation}

\noindent
where $\Lambda(n)$ is the Mangoldt function equal to $\ln(p)$ if $n=p^k$ for some prim $p$ and integer $k \geq 1$ and $0$ otherwise.

\vspace{2mm}
\noindent
Say we have a function $f$ holomorphic in the neighborhood of a point $a$. Suppose there is a zero in $a$, hence we can write $f(z)= (z - a)^n \, g(z)$ where $g(a) \neq 0$. From the derivative of $f$, we get $\varepsilon \frac{f^{\prime}(a + \varepsilon)}{f(a + \varepsilon)} = n + \varepsilon \frac{g^{\prime}(a+ \varepsilon)}{g(a+\varepsilon)}$. Suppose $a$ is not a pole of $f$, hence we get $\frac{g^{\prime}(a)}{g(a)}$ is equal to a constant. By taking the limit of the real part of $\varepsilon \frac{f^{\prime}(a+\varepsilon)}{f(a+\varepsilon)}$ when $\varepsilon$ tends to zero, we get:

\begin{equation}
w(f,a) = \lim\limits_{\varepsilon \to 0} \, \Re \left[ \varepsilon \, \frac{f^{\prime}(a + \varepsilon)}{f(a+ \varepsilon)} \right] = n \,,
\end{equation}

\noindent
which is the multiplicity of the zero in $a$.

\vspace{2mm}
\noindent
Let us us set set $s = 1 + \varepsilon + i \, \beta$ where $\beta$ is a real number and $\varepsilon > 0$. We get:

\begin{equation}
\Re \left[ \varepsilon \frac{\zeta^{\prime}(s)}{\zeta(s)} \right] = - \varepsilon \, \sum_{n=2}^{\infty} \Lambda(n) \, n^{-1 -\varepsilon} \, \cos \left( \beta \, \ln(n) \right) \,.
\end{equation}

\noindent
The proof uses the Mertens' trick, which is based on the below inequality:

\begin{equation}
\begin{split}
0 \leq \left( 1 + \cos \theta \right)^2 & = 1 + 2 \, \cos \theta + \cos^2 \theta \\
& = 1 + 2 \, \cos \theta + \frac{1 + \cos(2 \theta)}{2} \\
& = \frac{3}{2} + 2 \, \cos \theta + \frac{1}{2} \cos(2 \theta) \,.
\end{split}
\end{equation}

\noindent
Hence, we get:

\begin{equation}
3 + 4 \, \cos \theta + \cos(2\theta) \geq 0 \,.
\end{equation}

\noindent
From (40) and (41), we get:

\begin{equation}
3 \, \Re \left[\varepsilon \frac{\zeta^{\prime}}{\zeta}(1 + \varepsilon) \right] + 4 \Re \left[ \varepsilon \frac{\zeta^{\prime}}{\zeta} \left( 1 + \varepsilon + 2 \, i \, \beta \right) \right] + \Re \left[ \varepsilon \frac{\zeta^{\prime}}{\zeta} \left( 1 + \varepsilon + 2 \, i \beta \right)\right] \leq 0 \,,
\end{equation}

\noindent
where $\varepsilon > 0$. By taking the limit of (43) when $\varepsilon$ tends to zero from the right $(\varepsilon > 0)$, we get:

\begin{equation}
3 \, w \left(\zeta, 1 \right) + 4 \, w \left(\zeta, 1 + i \, \beta \right) + w \left(\zeta, 1 + 2 \, i \beta \right) \leq 0 \,.
\end{equation}

\noindent
Since $\zeta(s)$ only has a simple pole in $s = 1$, when $\beta \neq 0$, we have $w \left( \zeta, 1 + 2 \, i \beta\right) \geq 0$ and $w \left(\zeta,1\right) = -1$. Thus, we get:

\begin{equation}
\begin{split}
0 \leq w \left(\zeta, 1 + i \, \beta \right) & \leq - \frac{3}{4} \, w \left(\zeta, 1 \right) \\
& \leq \frac{3}{4} \,. 
\end{split}
\end{equation}

\noindent
As the Riemann zeta function is holomorphic, the multiplicity of a zero in $a$ denoted $w\left( \zeta, a \right)$ must be an integer. According to (45), $\forall \beta \neq 0, w(\zeta, 1+ i \, \beta)=0$. Thus, we can say that the Riemann zeta function has no zeros on the line $\Re(s)=1$.

\noindent
The fact that $\zeta(s)$ has no zeros on the line $\Re(s) = 0$ is obtained by reflecting the line $\Re(s)=1$ with the Riemann zeta functional. To show that $\zeta(0)=-1/2$, we use the Riemann zeta functional (3), which is expressed as:

\begin{equation}
\zeta(s) = \pi^{s-1} 2^s \, \sin \left( \frac{s \pi}{2}\right) \, \Gamma(1-s) \, \zeta(1-s) \,.
\end{equation}

\noindent
The Gamma functional equation, obtained using integration by part on integral form of Gamma function, see formula 6.1.15 p. 256 in \cite{Abramowitz} is expressed as follows:

\begin{equation}
s \, \Gamma(s) = \Gamma(s+1) \,.
\end{equation}

\noindent
From the multiplication of (46) with $(1-s)$, we get $(1-s) \zeta(s) = \pi^{s-1} \, 2^s \, \sin\left( \frac{s \pi}{2} \right) \, (1 - s) \, \Gamma(1-s)\, \zeta(1-s)$. From (47) we get $(1-s) \, \Gamma(1-s) = \Gamma(2-s)$. As $\zeta(s)$ has a simple pole in $s=1$, we get $\lim_{s \to 1} \left(1-s \right)\Gamma(s) = -1$. Also we have $\lim_{s \to 1} \pi^{s-1} \, 2^s \sin\left( \frac{s \pi}{2} \right) \, \Gamma(2-s) \, \zeta(1-s)= 2\, \zeta(0)$. Hence $\zeta(0)=-1/2$.

\section{Zeros in the critical strip $\Re(s) \in \, ]0,1[$}

The Riemann hypothesis states that all non-trivial zeros, those in the critical strip $\Re(s) \in \, ]0,1[$ lie on the critical line $\Re(s)=1/2$. By reflection of the zeros of the Riemann zeta function with respect to the critical line $\Re(s)=1/2$ (\textit{prop. no 8}), one has to show there are no zeros on either sides of the critical strip, either to the left or the right of the critical line $\Re(s)=1/2$, to say the Riemann hypothesis is true (an approach referred to as the one-sided RH test).

\bibliographystyle{spmpsci}      
\bibliography{Zeta}   

\begin{thebibliography}{1}
\providecommand{\url}[1]{{#1}}
\providecommand{\urlprefix}{URL }
\expandafter\ifx\csname urlstyle\endcsname\relax
  \providecommand{\doi}[1]{DOI~\discretionary{}{}{}#1}\else
  \providecommand{\doi}{DOI~\discretionary{}{}{}\begingroup
  \urlstyle{rm}\Url}\fi

\bibitem{Abramowitz}
Abramowitz, M., Stegun, I.: {H}andbook of {M}athematical {F}unctions with {F}ormulas, {G}raphs, and {M}athematical {T}ables.
\newblock National Bureau of Standards Applied Mathematics Series, Series \textbf{55} (1964)

\bibitem{Edwards}
Edwards, H.M.: {R}iemann's Zeta Function.
\newblock Dover Publications Inc. (2003)

\bibitem{Grigorieva}
Grigorieva, E.: {M}ethods of Solving Sequence and Series Problems,
\newblock Birkh{\"a}user (2016)

\bibitem{Hadamard}
Hadamard, J.: {S}ur la distribution des z{\'e}ros de la fonction $\zeta(s)$ et ses cons{\'e}quences arithm{\'e}tiques,
\newblock Bull. Soc. Math. France \textbf{24} pp. 199--220 (1896)

\bibitem{Hardy}
Hardy, G. H.: {R}amanujan twelve lectures on subjects suggested by his life and work.
\newblock AMS Chelsea publishing (1978); Hardy, Godfrey Harold (1877-1947)

\bibitem{HardyLittlewood}
Hardy, G.H., Littlewood, J.E.: {C}ontributions to the theory of the {R}iemann zeta function and the theory of the distribution of primes,
\newblock Acta Math. \textbf{41} pp. 119--196 (1916)

\bibitem{Ingham}
Ingham, A.E.: The Distribution of Prime Numbers.
\newblock Cambridge University Press, 1st published 1932 (1990)

\bibitem{Ivic}
Ivic, A.: The Riemann zeta-function: theory and applications.
\newblock Dover Publications, Inc (2003)

\bibitem{Krantz}
Krantz, S.G.: Foundations of Analysis.
\newblock Dover Publications, CRC Press (2015)

\bibitem{Poussin}
de la Vall{\'e}e Poussin, C.J.: {R}echerches analytiques sur la th{\'e}orie des nombres premiers.
\newblock Ann. Soc. Sci. Bruxelles \textbf{20}, 183--256 (1896)

\bibitem{RamanujanCollPapers}
Ramanujan, S.: Collected Papers of Srinivasa Ramanujan.
\newblock Cambridge University Press, edited by G.H. Hardy, P. V. Seshu Aiyar and B. M. Wilson (1927)

\bibitem{Ranjan}
Ranjan, R.: {T}he discovery of the {S}eries formula for $\pi$ by {L}eibniz, {G}regory and {N}ilakantha.
\newblock Mathematical Magazine \textbf{63} (5), 291--306 (1990)

\bibitem{Riemann1859}
Riemann, G.F.B.: {\"U}ber die {A}nzahl der {P}rimzahlen unter einer gegebenen {G}r{\"o}sse.
\newblock Monatsber. K{\"o}nigl. Preuss. Akad. Wiss. Berlin \textbf{Nov.},
  671--680 (1859)

\bibitem{Sondow}
Sondow, J.: {R}amanujan Primes and {B}ertrand's Postulate.
\newblock Amer. Math. Monthly \textbf{116},
  630--635 (2009)

\bibitem{SpencerandGraham}
Spencer, J., Graham, R.: {T}he elementary Proof of the Prime Number Theorem.
\newblock Monatsber. The Mathematical Intelligencer, Springer. \textbf{31},
  18--23 (2009)

\bibitem{Titchmarsh}
Titchmarsh, E.C., Heath-Brown, D.R.: The Theory of the Riemann Zeta-function.
\newblock Oxford University Press, 2nd ed. reprinted in 2007 (1986)

\bibitem{Veen}
van der Veen, R., van des Craat, J.: {T}he {R}iemann hypothesis.
\newblock MAA Press (2015)

\bibitem{Mangoldt}
von Mangoldt, H.: {Z}u {R}iemanns Abhandlung, {\"U}ber die {A}nzahl der {P}rimzahlen unter einer gegebenen {G}r{\"o}sse.
\newblock Journal f{\"u}r die reine und angewandte Mathematik \textbf{114},
  255--305 (1859)

\bibitem{Selberg}
Selberg, A.: {A}n elementary proof of the prime-number theorem.
\newblock {A}nnals of Mathematics, Second Series \textbf{50},
  305--313 (1949)

\end{thebibliography}

\end{document}